\documentclass[a4paper,11pt]{article}
\usepackage{stmaryrd}
\usepackage{graphicx}
\usepackage{amsfonts}
\usepackage{amsmath}
\usepackage{amssymb}
\usepackage{indentfirst}
\usepackage{titlesec}
\usepackage{caption2}
\usepackage{color}

\def\barr{\begin{array}}
\def\earr{\end{array}}
\def\bali{\begin{aligned}}
\def\eali{\end{aligned}}
\def\bearr{\begin{eqnarray}}
\def\eearr{\end{eqnarray}}

\providecommand{\play}{\displaystyle}

\providecommand{\li}{\limits}

\providecommand{\pt}{\partial}

\providecommand{\ra}{\rightarrow}

\providecommand{\da}{\downarrow}

\providecommand{\Prob}{\mathbf P}

\providecommand{\E}{\mathbf E}

\providecommand{\Var}{\mathbf {Var}}

\providecommand{\al}{\alpha}

\providecommand{\gm}{\gamma}

\providecommand{\dt}{\delta}

\providecommand{\Dt}{\Delta}

\providecommand{\ve}{\varepsilon}

\providecommand{\kp}{\kappa}

\providecommand{\sm}{\sigma}

\providecommand{\om}{\omega}

\providecommand{\Om}{\Omega}

\providecommand{\R}{\mathbb R}

\providecommand{\T}{\mathbb T}

\providecommand{\Z}{\mathbb Z}

\providecommand{\cE}{\mathcal E}

\providecommand{\cF}{\mathcal F}

\providecommand{\cG}{\mathcal G}

\providecommand{\cI}{\mathcal I}

\providecommand{\cR}{\mathcal R}

\providecommand{\grad}{\nabla}

\begin{document}

\title{It\^{o}'s formula, the stochastic exponential and change of measure on general time scales.}
\author{Wenqing Hu. \thanks{Department of Mathematics and Statistics, Missouri University of Science and Technology
(formerly University of Missouri, Rolla). Email: \texttt{huwen@mst.edu}}}

\date{
}

\maketitle

\begin{abstract}
We provide an It\^{o}'s formula for stochastic dynamical equation on general time scales. Based on this
It\^{o}'s formula we give a closed form expression for stochastic exponential on general time scales.
We then demonstrate a Girsanov's change of measure formula in the case of general time scales. Our result
is being applied to a Brownian motion on the quantum time scale ($q$--time scale).
\end{abstract}



\section{Introduction}

The theory of dynamical equation on time scales (\cite{[Bohner Peterson book]}) has attracted many researches recently.
In particular, attempts of extension to stochastic dynamical equations and stochastic analysis on general time scales
have been made in several previous works (\cite{[SDtE Bohner]}, \cite{[SD Exponential Bohner-Sanyal]}, \cite{[SdtE Grow-Sanyal]},
\cite{[time scale BM]}, \cite{[first attempt time scale stochastic]}).
In the work \cite{[SD Exponential Bohner-Sanyal]} the authors mainly work with a discrete time scale;
in \cite{[SDtE Bohner]} the authors introduce an extension of a function and define the stochastic as well as deterministic integrals
as the usual integrals for the extended function; in \cite{[SdtE Grow-Sanyal]} the authors make use
of their results on the quadratic variation of a Brownian motion (\cite{[quadratic variation BM time scale]})
on time scales and, based on this, they define the stochastic integral via a generalized version of the
It\^{o} isometry; in \cite{[first attempt time scale stochastic]} the authors introduce the so called $\grad$--stochastic
integral via the backward jump operator and they also derive an It\^{o}'s formula based on this definition
of the stochastic integral. We notice that different previous works adopt different notions of the stochastic integral and there lacks
a uniform and coherent theory of a stochastic calculus on general time scales.

The purpose of the present article is to fill in this gap. We will be mainly working under the framework of \cite{[SDtE Bohner]}, in that
we define our stochastic integral using the definition given in \cite{[SDtE Bohner]}. We then present a general It\^{o}'s formula
for stochastic dynamical equations under the framework of \cite{[SDtE Bohner]}. Our It\^{o}'s formula works for general time scales,
and thus fills the gap left in \cite{[SD Exponential Bohner-Sanyal]}, which deals with only discrete time scales. By making use of the
It\^{o}'s formula we obtain a closed--form expression for the stochastic exponential on general times scales. We will then
 demonstrate a change of measure (Girsanov's) theorem for stochastic dynamical equation on time scales.

We would like to point out that our change of measure formula is different from the continuous process case in that
the density function is not given by the stochastic exponential, but rather is found by the fact that the process on the time scale can be extended
to a continuous process simply by linear extension.

It is also worth mentioning that our construction is different from
\cite{[Bhamidi et al]} in that we are working with the case that the time parameter of the process is running on a time scale,
 whereas in \cite{[Bhamidi et al]} and related works (e.g. \cite{[Short time asymptotic BM fractal]},
 \cite{[LDP for BM fractal]}, \cite{[Cont-Fournie]}) the authors are working
with the case that the state space of the process is a time scale.

We note that stochastic calculus on the so called $q$--Brownian motion has been considered in \cite{[Haven2009]}, \cite{[Haven2011]},
\cite{[BrycqBM]}. As an application, we will also work our It\^{o}'s formula for a Brownian motion
on the quantum time scale ($q$--time scale)
case at the last section of the paper.

The paper is organized as follows. In Section 2 we discuss some basic set--up for time scales calculus. In Section 3 we will briefly review
the results in \cite{[SDtE Bohner]} and define the stochastic integral and stochastic dynamical equation on time scales. In Section 4
we present and prove our It\^{o}'s formula. In Section 5 we discuss the formula for stochastic exponential. In Section 6 we prove the change of
measure (Girsanov's) formula. Finally in Section 7 we consider an example of Brownian motion on a quantum time scale.

\section{Set--up: Basics of time scales calculus.}

A \textit{time scale} $\T$ is an arbitrary nonempty closed subset of the real numbers $\R$, where we assume that
$\T$ has the topology that it inherits from the real numbers $\R$ with the standard topology.

We define the \textit{forward jump operator} by
$$\sm(t)=\inf\{s\in \T: s>t\} \text{ for all } t\in \T \text{ such that this set is non--empty} \ ,$$
and the \textit{backward jump operator} by
$$\rho(t)=\sup\{s\in \T: s<t\} \text{ for all } t\in \T \text{ such that this set is non--empty} \ .$$

Let $t\in \T$.
If $\sm(t)>t$, then $t$ is called \textit{right--scattered}.
If $\sm(t)=t$, then $t$ is called \textit{right dense}.
If $\rho(t)<t$, then $t$ is called \textit{left--scattered}.
If $\rho(t)=t$, then $t$ is called \textit{left--dense}.
Moreover, the sets $\T^\kp$ and $\T_\kp$ are derived from $\T$ as follows: If $\T$
has a left--scattered maximum, then $\T^\kp$ is the set $\T$ without that left--scattered maximum;
otherwise, $\T^\kp=\T$. If $\T$ has a right--scattered minimum, then $\T_\kp$ is the set $\T$
without that right--scattered minimum; otherwise, $\T_\kp=\T$.
The \textit{graininess function} is defined by $\mu(t)=\sm(t)-t$ for all $t\in \T^\kp$.

Notice that since $\T$ is closed, for any $t\in \T$, the points $\sm(t)$ and $\rho(t)$ are belonging to $\T$.

For a set $A\subset \R$ we denote the set $A_\T=A\cap \T$.

\

Given a time scale $\T$ and a function $f: \T\ra \R$,
the \textit{delta} (or \textit{Hilger}) \textit{derivative} $f^\Dt(t)$ of $f$ at $t\in \T$
is defined as follows (\cite[Definition 1.10]{[Bohner Peterson book]}).

\

\textbf{Definition 1.}
\textit{Assume $f: \T\ra \R$ is a function and let $t\in \T^\kp$. Then we define $f^\Dt(t)$ to be the number (provided that it exists)
with the property that given any $\ve>0$, there is a neighborhood $U$ of $t$ (i.e., $U=(t-\dt, t+\dt)\cap \T$ for some $\dt>0$) such that}
$$\left|[f(\sm(t))-f(s)]-f^\Dt(t)[\sm(t)-s]\right|\leq \ve |\sm(t)-s| \text{ \textit{for all} } s\in U \ .$$

\

The delta derivative is characterized by the following Theorem, which is \cite[Theorem 1.16]{[Bohner Peterson book]}.

\

\textbf{Theorem 1.} \textit{Assume that $f: \T\ra \R$ is a function and let $t\in \T^\kp$. Then we have the following:}

(i) \textit{If $f$ is differentiable at $t$, then $f$ is continuous at $t$.}

(ii) \textit{If $f$ is continuous at $t$ and $t$ is right--scattered, then $f$ is differentiable at $t$ with}
$$f^\Dt(t)=\dfrac{f(\sm(t))-f(t)}{\sm(t)-t} \ .$$

(iii) \textit{If $t$ is right--dense, then $f$ is differentiable at $t$ if and only if the limit }
$$\lim\li_{s\ra t}\dfrac{f(t)-f(s)}{t-s}$$
\textit{exists as a finite number. In this case} $$f^\Dt(t)=\lim\li_{s\ra t}\dfrac{f(t)-f(s)}{t-s} \ .$$

(iv) \textit{If $f$ is differentiable at $t$, then} $$f(\sm(t))=f(t)+\mu(t)f^\Dt(t) \ .$$

\

\section{Stochastic integrals and stochastic differential equations on time scales.}

We will adopt the definitions introduced in \cite{[SDtE Bohner]}
as our definition of a Brownian motion and It\^{o}'s stochastic integral on time scales. In the next section we will derive
an It\^{o}'s formula corresponding to the stochastic integral defined in such a way.

\

\textbf{Definition 2.} \textit{A Brownian motion indexed by a time scale $\T\subset \R$ is an adapted stochastic process
$\{W_t\}_{t\in \T\cup \{0\}}$ on a filtered probability space $(\Om, \cF_t, \Prob)$ such that}

(1) \textit{$\Prob(W_0=0)=1$};

(2) \textit{If $s<t$ and $s,t\in \T$ then the increment $W_t-W_s$ is independent of $\cF_s$
and is normally distributed with mean $0$ and variance $t-s$};

(3) \textit{The process $W_t$ is almost surely continuous on $\T$.}

\

Note that the property (3) is proved in the work \cite{[time scale BM]}.

For a random function $f: [0,\infty)_{\T}\times \Om\ra \R$ we define the \textit{extension} $\widetilde{f}: [0,\infty)\times \Om\ra \R$
by

\begin{equation}
\label{Extension}\widetilde{f}(t,\om)=f(\sup[0,t]_{\T},\om)
\end{equation}
for all $t\in [0,\infty)$.

We shall make use of the definitions given in \cite{[SDtE Bohner]} for the classical Lebesgue and
Riemann integral. For any random function $f: [0,\infty)_{\T}\times \Om \ra \R$ and $T<\infty$
we define its $\Dt$--Riemann (Lebesgue) integral
as $$\int_0^T f(t,\om)\Dt t=\int_0^T \widetilde{f}(t,\om)dt \ ,$$ where the integral on the right hand side
of the above equation is interpreted as a standard Riemann (Lebesgue) integral. In a similar way, the work
\cite{[SDtE Bohner]} defines a stochastic integral for an $L^2([0,T]_\T)$--progressively measurable random function $f(t,\om)$
as $$\int_0^T f(t,\om)\Dt W_t=\int_0^T \widetilde{f}(t,\om)dW_t \ ,$$
where again the right hand side of the above equation is interpreted as a standard It\^{o}'s stochastic integral. Note that
the way (1) that we define the extension guarantees that the function $\widetilde{f}(t,\om)$ is progressively measurable.

In \cite{[SDtE Bohner]} the authors then defined the solution of the $\Dt$--stochastic differential equation indicated
by the notation

\begin{equation}\label{SDeltaE}
\Dt X_t=b(t,X)\Dt t+\sm(t,X)\Dt W_t \ ,
\end{equation}
as the process $\{X_t\}_{t\in [0,T]_\T}$ such that

$$X_{t_2}-X_{t_1}=\int_{t_1}^{t_2}b(t,X_t)\Dt t+\int_{t_1}^{t_2}\sm(t, X_t)\Dt W_t \ ,$$
with the deterministic and stochastic integrals on the right--hand side of the above equality
interpreted as was just mentioned. Under the condition of continuity in the $t$--variable
and uniform Lipschitz continuity in the $x$--variable of the functions
$b(t,x)$ and $\sm(t,x)$, together with being no worse than linear growth in $x$--variable, existence and pathwise uniqueness
of strong solution to \eqref{SDeltaE} are proved in \cite{[SDtE Bohner]}.

\

\section{It\^{o}'s formula for stochastic integrals on time scales.}

We will make use of the following fact that is simple to prove.

\

\textbf{Proposition 1.} \textit{The set of all left--scattered or right--scattered points of $\T$ is at most countable.}

\

\textbf{Proof.} If $x\in \T$ is a right--scattered point, then $I_x=(x,\sm(x))$ is an open interval such that $I_x\cap \T=\emptyset$.
Similarly, if $x\in \T$ is a left--scattered point, then $I_x=(\rho(x), x)$ is an open interval such that $I_x\cap \T=\emptyset$.
Suppose $x<y$ and $x,y\in \T$. We then distinguish four different cases.

Case 1: both $x$ and $y$ are right--scattered. We argue that in this case we have $I_x\cap I_y=\emptyset$.
Suppose this is not the case, then we must have $\sm(x)>y$. But we see that $\sm(x)=\inf\{s>x: s\in \T\}$ and
$y\in \T$. So we must have $\sm(x)\leq y$. We arrive at a contradiction;

Case 2: both $x$ and $y$ are left--scattered. This case is similar to Case 1 and we conclude that $I_x\cap I_y=\emptyset$;

Case 3: $x$ is left--scattered, $y$ is right--scattered. In this case we see that $I_x=(\rho(x), x)$ and $I_y=(y, \sm(y))$,
as well as $x<y$. This implies that $I_x\cap I_y=\emptyset$;

Case 4: $x$ is right--scattered, $y$ is left--scattered. In this case $I_x=(x, \sm(x))$ and $I_y=(\rho(y), y)$.
If $\sm(x)\leq\rho(y)$, then $I_x\cap I_y=\emptyset$. If $\sm(x)>\rho(y)$, then we see that $(x,y)=I_x\cup I_y$ so that
$(x,y)\cap \T=\emptyset$. That implies further $\sm(x)=y$ and $\rho(y)=x$, i.e., $I_x=I_y$.

Thus we see that for all points $x\in \T$ being left or right--scattered, the set of
all open intervals of the form $I_x$ are disjoint subsets of $\R$. Henceforth there are at most countably many such intervals.
Each such interval corresponds to one or two endpoints in $\T$ that are either left or right--scattered. Thus the total number
of left or right--scattered points in $\T$ are at most countably many. $\square$

\

Let $C$ be the (at most) countable set of all left--scattered or right--scattered points of $\T$.
As we have already seen in the proof of the previous Proposition, the set $C$ corresponds to at most countably many
open intervals $\cI=\{I_1, I_2,...\}$ such that (1) for any $k\neq l$, $I_k\cap I_l=\emptyset$; (2) either the left endpoint
or right endpoint or both endpoints of any of the $I_k$'s are in $\T$, and are left or right--scattered; (3) $I_k\cap \T=\emptyset$
for any $k=1,2,...$; (4) any point in $C$ is a left or right--endpoint of one of the $I_k$'s.

We will denote $I_k=(s_{I_k}^-, s_{I_k}^+)$. Since for any $x\in \T$, the points $\sm(x)$ and $\rho(x)$
are in $\T$, we further infer that for any such interval $I_k$, we have $s_{I_k}^-$ and $s_{I_k}^+$ are in $\T$,
so that $s_{I_k}^-$ is right--scattered and $s_{I_k}^+$ is left--scattered.

We then establish the following It\^{o}'s formula.

For any two points $t_1, t_2\in \T$, $t_1\leq t_2$, and any open interval $I_k\in \cI$, such that
$I_k\cap [t_1, t_2]\neq \emptyset$, we have $I_k\subset (t_1,t_2)$. This is because if not the case, then $t_1$
or $t_2$ will belong to $I_k$, contradictory to the fact that $I_k\cap \T=\emptyset$. We conclude that

$$\{I_k \in \cI: I_k\cap [t_1,t_2]\neq \emptyset\}=\{I_k\in \cI: I_k\subset (t_1, t_2)\} \ .$$

Let us consider a function $f(t,x): \T\times \R\ra \R$. Let $f^\Dt(t,x)$, $f^{\Dt^2}(t,x)$
be the first and second order delta (Hilger) derivatives of $f$
with respect to time variable $t$ at $(t,x)$ and
let $\dfrac{\pt f}{\pt x}(t,x)$ and $\dfrac{\pt^2 f}{\pt x^2}(t,x)$ be the first and second order partial derivatives of $f$
with respect to space variable $x$ at $(t,x)$.

\

\textbf{Theorem 2.} (It\^{o}'s formula) \textit{Let any function $f: \T\times \R\ra \R$ be such that
$f^\Dt(t,x)$, $f^{\Dt^2}(t,x)$, $\dfrac{\pt f}{\pt x}(t,x)$, $\dfrac{\pt^2 f}{\pt x^2}(t,x)$,
$\dfrac{\pt f^\Dt}{\pt x}(t,x)$ and $\dfrac{\pt^2 f^\Dt}{\pt x^2}(t,x)$
are continuous on $\T\times \R$.
Set any $t_1\leq t_2$, $t_1, t_2\in [0, \infty)_\T$, then we have}

\begin{equation}\label{ItoFormula}\begin{array}{ll}
&f(t_2,W_{t_2})-f(t_1, W_{t_1})
\\
=&\play{\int_{t_1}^{t_2}f^\Dt(s, W_s)\Dt s+
\int_{t_1}^{t_2}\dfrac{\pt f}{\pt x}(s,W_s)\Dt W_s
+\dfrac{1}{2}\int_{t_1}^{t_2}\dfrac{\pt^2 f}{\pt x^2}(s, W_s)\Dt s}
\\
&\qquad +\play{\sum\li_{I_k\in \cI, I_k\subset (t_1,t_2)} \left[f(s_{I_k}^+,W_{s_{I_k}^+})-f(s_{I_k}^+, W_{s_{I_k}^-})
\right.}
\\
&\play{\qquad \qquad \qquad \qquad \qquad \left.-\dfrac{\pt f}{\pt x}(s_{I_k}^-,W_{s_{I_k}^-})(W_{s_{I_k}^+}-W_{s_{I_k}^-})
-\dfrac{1}{2}\dfrac{\pt^2 f}{\pt x^2}(s_{I_k}^-,W_{s_{I_k}^-})(s_{I_k}^+-s_{I_k}^-)\right] \ .}
\end{array}\end{equation}

\

\textbf{Proof.} We will make use of the following classical version (Peano form) of Taylor's theorem: for any function
$f: \T\times \R\ra \R$ such that $\dfrac{\pt f}{\pt x}(t,x)$ and $\dfrac{\pt^2 f}{\pt x^2}(t,x)$ are continuous on $\T\times \R$,
and any $s\in \T$ and $x_1, x_2\in \R$ we have

\begin{equation}\label{ClassicalTaylorExpansion}
f(s, x_2)-f(s, x_1)=\dfrac{\pt f}{\pt x}(s,x_1)(x_2-x_1)
+\dfrac{1}{2}\dfrac{\pt^2 f}{\pt x^2}(s,x_1)(x_2-x_1)^2+R^f_\text{C}(s; x_1,x_2) \ ,
\end{equation}
where $$|R^f_\text{C}(s; x_1, x_2)|\leq r(|x_2-x_1|)(x_2-x_1)^2 \ ,$$ and $r: \R_+\ra \R_+$
is an increasing function with $\lim\li_{u\da 0}r(u)=0$.

We will also make use of the time scale Taylor formula (see \cite[Theorem 1.113]{[Bohner Peterson book]} as well as
\cite{[AgarwalBohnerBasicTimeScaleCalc]}) applied to $f(t,x)$ up to first order in $t$:
for any $s_2>s_1$ and $s_1,s_2\in \T$ we have

\begin{equation}\label{TimeScaleTaylorExpansion}
f(s_2, x)-f(s_1, x)= f^{\Dt}(s_1,x)(s_2-s_1)+R^f_{\text{TS}}(x; s_1,s_2) \ ,
\end{equation}
where $$|R^f_{\text{TS}}(x; s_1,s_2)|=\left|\int_{s_1}^{\rho(s_2)}(s_2-\sm(s))f^{\Dt^2}(s)\Dt s\right|\leq r(|s_2-s_1|)|s_2-s_1|$$
with $r(\bullet)$ as before.

Combining \eqref{ClassicalTaylorExpansion} and \eqref{TimeScaleTaylorExpansion} we see that we have

\begin{equation}\label{TaylorExpansion}
\begin{array}{ll}
& f(s_2, x_2)-f(s_1, x_1)
\\
=& [f(s_2,x_2)-f(s_1,x_2)]+[f(s_1,x_2)-f(s_1,x_1)]
\\
=& f^{\Dt}(s_1,x_2)(s_2-s_1)+\dfrac{\pt f}{\pt x}(s_1,x_1)(x_2-x_1)
+\dfrac{1}{2}\dfrac{\pt^2 f}{\pt x^2}(s_1,x_1)(x_2-x_1)^2
+R^f_{\text{TS}}(x_2;s_1,s_2)+R^f_\text{C}(s_1; x_1,x_2)
\\
=& \left[f^{\Dt}(s_1,x_1)+\dfrac{\pt f^\Dt}{\pt x}(s_1,x_1)(x_2-x_1)+\dfrac{1}{2}\dfrac{\pt^2 f^\Dt}{\pt x^2}(s_1,x_1)(x_2-x_1)^2
+R_C^{f^\Dt}(s_1; x_1,x_2)\right](s_2-s_1)
\\
& \qquad +\dfrac{\pt f}{\pt x}(s_1,x_1)(x_2-x_1)+\dfrac{1}{2}\dfrac{\pt^2 f}{\pt x^2}(s_1,x_1)(x_2-x_1)^2
+R^f_{\text{TS}}(x_2;s_1,s_2)+R^f_\text{C}(s_1; x_1,x_2)
\\
=& f^{\Dt}(s_1,x_1)(s_2-s_1)+\dfrac{\pt f}{\pt x}(s_1,x_1)(x_2-x_1)
+\dfrac{1}{2}\dfrac{\pt^2 f}{\pt x^2}(s_1,x_1)(x_2-x_1)^2
+R(s_1,s_2;x_1,x_2) \ ,
\end{array}
\end{equation}
with $$|R(s_1,s_2;x_1,x_2)|\leq r(|s_2-s_1|)|s_2-s_1|+r(|x_2-x_1|)(x_2-x_1)^2$$
for another function $r: \R_+\ra \R_+$ increasing with $\lim\li_{u\da 0}r(u)=0$.

Consider a partition $\pi^{(n)}$: $t_1=s_0<s_1<...<s_n=t_2$, such that (1) each $s_i\in \T$;
(2) $\max\li_{i}(\rho(s_i)-s_{i-1})\leq \dfrac{1}{2^n}$
for $i=1,2,...,n$. Notice that by definition $\rho(s_i)=\sup\{s<s_i: s\in \T\}$, so that we can always find $s_{i-1}\in \T$
so that $\rho(s_i)-s_{i-1}$ is sufficiently small.

Let the sets $C$ and $\cI$ be defined as before.
Let us fix a partition $\pi^{(n)}$, and consider a classification of its corresponding intervals $(s_{i-1}, s_i) \ , \ i=1,2,...,n$.
We will classify all intervals $(s_{i-1}, s_i)$ such that for all $I_k\in \cI$ we have $I_k\cap (s_{i-1}, s_i)=\emptyset$
as class $(a)$; and we classify
all intervals $(s_{i-1}, s_i)$ such that there exist some $I_k\in \cI$ with
$(s_{i-1}, s_i)\cap I_k \neq \emptyset$ as class $(b)$. For an interval $(s_{i-1}, s_i)$ in class $(a)$,
since for all $I_k\in \cI$ we have $I_k\cap (s_{i-1}, s_i)=\emptyset$, we see that $\rho(s_i)=s_i$. Because otherwise $(\rho(s_i), s_i)$
will be one of the $I_k$'s. Thus in this case we have $s_i-s_{i-1}<\dfrac{1}{2^n}$.
For an interval $(s_{i-1}, s_i)$ in class $(b)$,
since both $s_{i-1}$ and $s_i$
are in $\T$, we see that we have in fact $I_k\subseteq (s_{i-1}, s_i)$. In this case either $I_k=(s_{i-1}, s_i)$,
or $I_k\neq (s_{i-1}, s_i)$. If the latter happens,
then $(\rho(s_i), s_i)\in \cI$ is one of the $I_k$'s and $\rho(s_i)-s_{i-1}<\dfrac{1}{2^n}$.
We also see from the above analysis that all $I_k$'s are contained in intervals $(s_{i-1}, s_i)$ that belong to
class $(b)$. On the other hand, each interval $(s_{i-1}, s_i)$ is either entirely one of the $I_k$'s, or
it contains an interval $(\rho(s_i), s_i)$ that is one of the $I_k$'s. For the latter case, i.e., when
$s_{i-1}<\rho(s_i)<s_i$,
the set of intervals of the form $(s_{i-1}, \rho(s_i))$ are disjoint open intervals such that
\begin{equation}\label{EstimateOfRemainderIntervals}
\sum\li_{(s_{i-1}, s_i)\in (b), \, s_{i-1}<\rho(s_i)<s_i}(\rho(s_i)-s_{i-1})<\dfrac{n}{2^n} \ .
\end{equation}

Now we have

\begin{equation}\label{ItoFormulaExpansion}
\begin{array}{l}
f(t_2, W_{t_2})-f(t_1, W_{t_1})
\\
\play{=\sum\li_{i=1}^n[f(s_i, W_{s_i})-f(s_{i-1}, W_{s_{i-1}})]}
\\
\play{=\sum\li_{(a)}[f(s_i, W_{s_i})-f(s_{i-1}, W_{s_{i-1}})]
+\sum\li_{(b)}[f(s_i, W_{s_i})-f(s_{i-1}, W_{s_{i-1}})]}
\\
\play{=(I)+(II) \ .}
\end{array}
\end{equation}

We apply \eqref{TaylorExpansion} term by term in part $(I)$ of \eqref{ItoFormulaExpansion}, and we get

\begin{equation}
\begin{array}{l}\label{TaylorExpansionOfII}
\play{(I)}
\\
\play{=\sum\li_{(a)}[f(s_i, W_{s_i})-f(s_{i-1}, W_{s_{i-1}})]}
\\
\play{=\sum\li_{(a)}\left[f^\Dt(s_{i-1}, W_{s_{i-1}})(s_i-s_{i-1})
+\dfrac{\pt f}{\pt x}(s_{i-1}, W_{s_{i-1}})(W_{s_i}-W_{s_{i-1}})\right.}
\\
\play{\left. \ \ \ \ \ \ \ \ \ \ \ \ \ \ \ \ \ \ \ \ \ \ \ \ \ \ \ \ \ \ \ \ \ \ \ \ \
+\dfrac{1}{2}\dfrac{\pt^2 f}{\pt x^2}(s_{i-1}, W_{s_{i-1}})(W_{s_i}-W_{s_{i-1}})^2+R(s_{i-1}, s_i; W_{s_{i-1}}, W_{s_i})\right]}
\\
\play{=\sum\li_{i=1}^n \left[f^\Dt(s_{i-1}, W_{s_{i-1}})(s_i-s_{i-1})
+\dfrac{\pt f}{\pt x}(s_{i-1}, W_{s_{i-1}})(W_{s_i}-W_{s_{i-1}})\right]}
\\
\play{\ \ \ \ \ \ +\left(\sum\li_{(a)}\dfrac{1}{2}\dfrac{\pt^2 f}{\pt x^2}(s_{i-1}, W_{s_{i-1}})(W_{s_i}-W_{s_{i-1}})^2
                        +\sum\li_{(b)}\dfrac{1}{2}\dfrac{\pt^2 f}{\pt x^2}(s_{i-1}, W_{s_{i-1}})(s_i-s_{i-1})\right)}
\\
\play{\ \ \ \ \ \ -\sum\li_{(b)}\left[f^\Dt(s_{i-1}, W_{s_{i-1}})(s_i-s_{i-1})
+\dfrac{\pt f}{\pt x}(s_{i-1}, W_{s_{i-1}})(W_{s_i}-W_{s_{i-1}})+\dfrac{1}{2}\dfrac{\pt^2 f}{\pt x^2}(s_{i-1}, W_{s_{i-1}})(s_i-s_{i-1})\right]}
\\
\play{\ \ \ \ \ \ +\sum\li_{(a)}R(s_{i-1}, s_i; W_{s_{i-1}}, W_{s_i}) \ .}
\\
\play{= (III)_1+(III)_2+(IV)+(V) \ .}
\end{array}
\end{equation}

We have the following four convergence results:

\textit{Convergence Result} 1.1. By Lemma 1 \eqref{ConvergenceDeterministicDeltaIntegral} and \eqref{ConvergenceStochasticDeltaIntergal}
established below we have

\begin{equation}\label{ConvergenceResult1-1}
\Prob\left((III)_1\ra \play{\int_{t_1}^{t_2}f^\Dt(s, W_s)\Dt s+\int_{t_1}^{t_2}\dfrac{\pt f}{\pt x}(s, W_s)\Dt W_s}
\text{ as $n \ra \infty$}\right)=1 \ .
\end{equation}

\textit{Convergence Result} 1.2. By Lemma 2 \eqref{ContinuousPartQuadraticVariation} and Lemma 1 \eqref{ConvergenceDeterministicDeltaIntegral}
established below we have

\begin{equation}\label{ConvergenceResult1-2}
\Prob\left((III)_2\ra \play{\int_{t_1}^{t_2}\dfrac{1}{2}\dfrac{\pt^2 f}{\pt x^2}(s, W_s)\Dt s}
\text{ as $n \ra \infty$}\right)=1 \ .
\end{equation}

\textit{Convergence Result} 2. We have, with probability $1$, that
$$(V)=\sum\li_{(a)}R(s_{i-1}, s_i; W_{s_{i-1}}, W_{s_i})\ra 0$$
as $n \ra \infty$.

In fact, by the Kolmogorov--\v{C}entsov theorem proved in Theorem 3.1 of \cite{[time scale BM]}
we know that for almost all trajectory of $W_t$ on $\T$, for each fixed trajectory $W_t(\om)$,
there exist an $n_0=n_0(\om)$ such that for all $n\geq n_0$, for a partition $\pi^{(n)}$ with a classification
of its intervals $(s_{i-1}, s_i)$ into classes $(a)$ and $(b)$ as above, $\sup\li_{(a)}|W_{s_i}-W_{s_{i-1}}|\leq \dfrac{\dt}{2^{\gm n/5}}$
for some fixed $\dt>0$ and $\gm>0$.  From here we can estimate

$$\begin{array}{l}
\play{\E\sum\li_{(a)}R(s_{i-1}, s_i; W_{s_{i-1}}, W_{s_i})}
\\
\leq \play{\E\sum\li_{(a)}\left[r(s_i-s_{i-1})(s_i-s_{i-1})+r(|W_{s_i}-W_{s_{i-1}}|)(W_{s_i}-W_{s_{i-1}})^2\right]}
\\
\leq r\left(\dfrac{1}{2^n}\right)(t_2-t_1)+r\left(\dfrac{\dt}{2^{\gm n/5}}\right)(t_2-t_1) \ , \end{array}$$
i.e.,

\begin{equation}\label{ConvergenceResult2}
\Prob\left(\lim\li_{n \ra \infty}\sum\li_{(a)}R(s_{i-1}, s_i; W_{s_{i-1}}, W_{s_i})=0\right)=1 \ .
\end{equation}

\textit{Convergence Result} 3. Let

$$\begin{array}{ll}
(II)+(IV)=A_n=A_n(\om)&\play{=\sum\li_{(b)}\left[f(s_i, W_{s_i})-f(s_{i-1}, W_{s_{i-1}})-f^\Dt(s_{i-1}, W_{s_{i-1}})(s_i-s_{i-1})\right.}
\\
&\play{\left. \ \ \ \ \ \
-\dfrac{\pt f}{\pt x}(s_{i-1}, W_{s_{i-1}})(W_{s_i}-W_{s_{i-1}})-\dfrac{1}{2}\dfrac{\pt^2 f}{\pt x^2}(s_{i-1}, W_{s_{i-1}})(s_i-s_{i-1})\right]
 \, ,}
\end{array}$$
and
$$\begin{array}{ll}
B_n=B_n(\om)&=\play{\sum\li_{I_k\in \cI, I_k\subset (t_1,t_2)} \left[f(s_{I_k}^+,W_{s_{I_k}^+})-f(s_{I_k}^-, W_{s_{I_k}^-})
-f^\Dt(s_{I_k}^-, W_{s_{I_k}^-})(s_{I_k}^+-s_{I_k}^-)
\right.}
\\
&\play{\qquad \qquad \qquad \qquad \left.-\dfrac{\pt f}{\pt x}(s_{I_k}^-,W_{s_{I_k}^-})(W_{s_{I_k}^+}-W_{s_{I_k}^-})
-\dfrac{1}{2}\dfrac{\pt^2 f}{\pt x^2}(s_{I_k}^-,W_{s_{I_k}^-})(s_{I_k}^+-s_{I_k}^-)\right]} \ .
\end{array}$$

We claim that we have

\begin{equation}\label{ConvergenceResult3}
\Prob(|A_n(\om)-B_n(\om)|\ra 0 \text{ as } n \ra \infty)=1 \ .
\end{equation}
In fact, from the analysis that leads to the estimate $(5)$ we see that we can write $A_n-B_n$
as

$$\begin{array}{ll}
&A_n-B_n
\\
=&\play{
\sum\li_{(s_{i-1}, s_i)\in (b), \, s_{i-1}<\rho(s_i)<s_i}
\left[f(s_i, W_{s_i})-f(s_{i-1}, W_{s_{i-1}})-f^\Dt(s_{i-1}, W_{s_{i-1}})(s_i-s_{i-1})\right.}
\\
&\play{\left. \ \ \ \ \ \
-\dfrac{\pt f}{\pt x}(s_{i-1}, W_{s_{i-1}})(W_{s_i}-W_{s_{i-1}})-\dfrac{1}{2}\dfrac{\pt^2 f}{\pt x^2}(s_{i-1}, W_{s_{i-1}})(s_i-s_{i-1})\right]}
\\
&\play{-\sum\li_{(s_{i-1}, s_i) \in (b), \, s_{i-1}<\rho(s_i)<s_i}
\left[f(s_i, W_{s_i})-f(\rho(s_i), W_{\rho(s_i)})-f^\Dt(\rho(s_i), W_{\rho(s_i)})(s_i-\rho(s_i))\right.}
\\
&\play{\left. \ \ \ \ \ \
-\dfrac{\pt f}{\pt x}(\rho(s_i), W_{\rho(s_i)})(W_{s_i}-W_{\rho(s_i)})
-\dfrac{1}{2}\dfrac{\pt^2 f}{\pt x^2}(\rho(s_i), W_{\rho(s_i)})(s_i-\rho(s_i))\right]}
\\
&\play{-\sum\li_{I_k\in \cI, I_k\subset (s_{i-1}, \rho(s_i)) \text{ for some } (s_{i-1}, s_i)
\in (b) } \left[f(s_{I_k}^+,W_{s_{I_k}^+})-f(s_{I_k}^-, W_{s_{I_k}^-})
\right.}
\\
&\play{\ \ \ \ \ \ \left.-f^\Dt(s_{I_k}^-, W_{s_{I_k}^-})(s_{I_k}^+-s_{I_k}^-)-\dfrac{\pt f}{\pt x}(s_{I_k}^-,W_{s_{I_k}^-})(W_{s_{I_k}^+}-W_{s_{I_k}^-})
-\dfrac{1}{2}\dfrac{\pt^2 f}{\pt x^2}(s_{I_k}^-,W_{s_{I_k}^-})(s_{I_k}^+-s_{I_k}^-)\right]}
\\
=&(VI)_1+(VI)_2+(VI)_3+(VI)_4-(VII) \ .
\end{array}$$

Here

$$(VI)_1=\sum\li_{(s_{i-1}, s_i)\in (b), \, s_{i-1}<\rho(s_i)<s_i}
[f(\rho(s_i), W_{\rho(s_i)})-f(s_{i-1}, W_{s_{i-1}})]  \ , $$

$$\begin{array}{ll}
(VI)_2&=\sum\li_{(s_{i-1}, s_i)\in (b), \, s_{i-1}<\rho(s_i)<s_i}
\left[f^\Dt(\rho(s_i), W_{\rho(s_i)})(s_i-\rho(s_i))
-f^\Dt(s_{i-1}, W_{s_{i-1}})(s_i-s_{i-1})\right]
\\
&=\sum\li_{(s_{i-1}, s_i)\in (b), \, s_{i-1}<\rho(s_i)<s_i}
\left[\left(f^\Dt(\rho(s_i), W_{\rho(s_i)})-f^\Dt(s_{i-1}, W_{s_{i-1}})\right)(s_i-s_{i-1})
\right.
\\
& \ \ \ \ \ \ \ \ \ \ \ \ \ \ \ \ \ \ \ \ \ \ \ \ \ \ \ \ \ \ \ \ \ \ \ \ \ \ \ \ \ \ \ \ \ \
\left.-f^\Dt(\rho(s_i), W_{\rho(s_i)})(\rho(s_i)-s_{i-1})\right] \ ,
\end{array}$$

$$\begin{array}{ll}
(VI)_3&=\sum\li_{(s_{i-1}, s_i)\in (b), \, s_{i-1}<\rho(s_i)<s_i}
\left[f^\Dt(\rho(s_i), W_{\rho(s_i)})(W_{s_i}-W_{\rho(s_i)})
-f^\Dt(s_{i-1}, W_{s_{i-1}})(W_{s_i}-W_{s_{i-1}})\right]
\\
&=\sum\li_{(s_{i-1}, s_i)\in (b), \, s_{i-1}<\rho(s_i)<s_i}
\left[\left(f^\Dt(\rho(s_i), W_{\rho(s_i)})-f^\Dt(s_{i-1}, W_{s_{i-1}})\right)(W_{s_i}-W_{s_{i-1}})
\right.
\\
& \ \ \ \ \ \ \ \ \ \ \ \ \ \ \ \ \ \ \ \ \ \ \ \ \ \ \ \ \ \ \ \ \ \ \ \ \ \ \ \ \ \ \ \ \ \
\left.-f^\Dt(\rho(s_i), W_{\rho(s_i)})(W_{\rho(s_i)}-W_{s_{i-1}})\right] \ ,
\end{array}$$

$$\begin{array}{ll}
(VI)_4&=\sum\li_{(s_{i-1}, s_i)\in (b), \, s_{i-1}<\rho(s_i)<s_i}
\left[f^\Dt(\rho(s_i), W_{\rho(s_i)})(s_i-\rho(s_i))
-f^\Dt(s_{i-1}, W_{s_{i-1}})(s_i-s_{i-1})\right]
\\
&=\sum\li_{(s_{i-1}, s_i)\in (b), \, s_{i-1}<\rho(s_i)<s_i}
\left[f^\Dt(s_{i-1}, W_{s_{i-1}})\left(s_{i-1}-\rho(s_i))\right)
\right.
\\
& \ \ \ \ \ \ \ \ \ \ \ \ \ \ \ \ \ \ \ \ \ \ \ \ \ \ \ \ \ \ \ \ \ \ \ \
\left.+\left(f^\Dt(\rho(s_i), W_{\rho(s_i)})-f^\Dt(s_{i-1}, W_{s_{i-1}})\right)
(s_i-\rho(s_i))\right] \ ,
\end{array}$$

$$\begin{array}{l}
(VII)
\\
=\play{\sum\li_{I_k\in \cI, I_k\subset (s_{i-1}, \rho(s_i)) \text{ for some } (s_{i-1}, s_i)
\in (b) } \left[f(s_{I_k}^+,W_{s_{I_k}^+})-f(s_{I_k}^-, W_{s_{I_k}^-})
\right.}
\\
\play{\ \ \ \ \ \ \left.-f^\Dt(s_{I_k}^-, W_{s_{I_k}^-})(s_{I_k}^+-s_{I_k}^-)-\dfrac{\pt f}{\pt x}(s_{I_k}^-,W_{s_{I_k}^-})(W_{s_{I_k}^+}-W_{s_{I_k}^-})
-\dfrac{1}{2}\dfrac{\pt^2 f}{\pt x^2}(s_{I_k}^-,W_{s_{I_k}^-})(s_{I_k}^+-s_{I_k}^-)\right]} \ .
\end{array}$$

From \eqref{EstimateOfRemainderIntervals}, the Kolmogorov--\v{C}entsov theorem proved in Theorem 3.1 of \cite{[time scale BM]},
as well as the assumptions about function $f$, we see that
$$\Prob\left(|(VI)_1|+|(VI)_2|+|(VI)_3|+|(VI)_4|+|(VII)|\ra 0 \text{ as } n \ra \infty\right)=1 \ .$$
From here we immediately see the claim \eqref{ConvergenceResult3}.

Note that for any interval $I_k=(s_{I_k}^-, s_{I_k}^+)$
we have $f^\Dt(s_{I_k}^-, W_{s_{I_k}^-})=\dfrac{f(s_{I_k}^+, W_{s_{I_k}^-})-f(s_{I_k^-}, W_{s_{I_k}^-})}{s_{I_k}^+-s_{I_k}^-}$, therefore
we see that

\begin{equation}\label{SubtractionItoBn}
\begin{array}{ll}
B_n=B_n(\om)&=\play{\sum\li_{I_k\in \cI, I_k\subset (t_1,t_2)} \left[f(s_{I_k}^+,W_{s_{I_k}^+})-f(s_{I_k}^-, W_{s_{I_k}^-})
-[f(s_{I_k}^+, W_{s_{I_k}^-})-f(s_{I_k^-}, W_{s_{I_k}^-})]
\right.}
\\
&\play{\qquad \qquad \qquad \qquad \left.-\dfrac{\pt f}{\pt x}(s_{I_k}^-,W_{s_{I_k}^-})(W_{s_{I_k}^+}-W_{s_{I_k}^-})
-\dfrac{1}{2}\dfrac{\pt^2 f}{\pt x^2}(s_{I_k}^-,W_{s_{I_k}^-})(s_{I_k}^+-s_{I_k}^-)\right]}
\\
&=\play{\sum\li_{I_k\in \cI, I_k\subset (t_1,t_2)} \left[f(s_{I_k}^+,W_{s_{I_k}^+})-f(s_{I_k}^+, W_{s_{I_k}^-})
\right.}
\\
&\play{\qquad \qquad \qquad \qquad \left.-\dfrac{\pt f}{\pt x}(s_{I_k}^-,W_{s_{I_k}^-})(W_{s_{I_k}^+}-W_{s_{I_k}^-})
-\dfrac{1}{2}\dfrac{\pt^2 f}{\pt x^2}(s_{I_k}^-,W_{s_{I_k}^-})(s_{I_k}^+-s_{I_k}^-)\right]} \ .
\end{array}
\end{equation}

Combining the convergence results \eqref{ConvergenceResult1-1}, \eqref{ConvergenceResult1-2}, \eqref{ConvergenceResult2}, \eqref{ConvergenceResult3},
together with \eqref{ItoFormulaExpansion} and
\eqref{TaylorExpansionOfII}, \eqref{SubtractionItoBn}, we establish \eqref{ItoFormula}. $\square$

\

The next two lemmas are used in the above proof of It\^{o}'s formula, but they are also of independent interest.

\

\textbf{Lemma 1.} (Convergence of $\Dt$--deterministic and stochastic integrals.)

\textit{Given a time scale $\T$ and $t_1, t_2\in \T$, $t_1<t_2$; a probability space
$(\Om, \cF, \Prob)$; a Brownian motion $\{W_t\}_{t\in \T}$ on the time scale $\T$,
for any progressively measurable random function $f$ that is continuous on $[t_1, t_2]\cap \T$,
viewed as a $L^2([t_1, t_2]_\T)$--progressively
 measurable random function $f(t, \om)$ on $\T$, and the families of partitions $\pi^{(n)}: t_1=s_0<s_1<...<s_n=t_2$,
$s_0, s_1, ...,s_n\in \T$, $\max\li_{i=1,2,...,n}(\rho(s_i)-s_{i-1})<\dfrac{1}{2^n}$, we have
}
\begin{equation}\label{ConvergenceDeterministicDeltaIntegral}
\Prob\left(\lim\li_{n\ra \infty}\sum\li_{i=1}^n f(s_{i-1}, \om)(s_i-s_{i-1})=\play{\int_{t_1}^{t_2}f(s,\om)\Dt s}\right)=1 \  ,
\end{equation}

\begin{equation}\label{ConvergenceStochasticDeltaIntergal}
\Prob\left(\lim\li_{n\ra \infty}\sum\li_{i=1}^n f(s_{i-1}, \om)(W_{s_i}-W_{s_{i-1}})=\play{\int_{t_1}^{t_2}f(s,\om)\Dt W_s}\right)=1 \  .
\end{equation}

\

\textbf{Proof.} As we have seen in the proof of the It\^{o}'s formula, for a given partition
$\pi^{(n)}: t_1=s_0<s_1<...<s_n=t_2$, such that $s_i\in \T$ for $i=0,1,...,n$, and $\max\li_{i=1,2,...,n}(\rho(s_i)-s_{i-1})<\dfrac{1}{2^n}$,
we can classify all intervals of the form $(s_{i-1}, s_i)$ into two classes $(a)$ and $(b)$: class $(a)$
are those open intervals $(s_{i-1}, s_i)$ such that it does not contain any open intervals
$I_k\in \cI$; class $(b)$ are those open intervals $(s_{i-1}, s_i)$ such that it contains at least one open interval
$I_k\in \cI$, the latter of which has endpoints that are left or right scattered.

Let us form a family of partitions $\sm^{(n)}: t_1=r_0<r_1<...<r_m=t_2$, so that the partition
$\sm^{(n)}$ is the partition $\pi^{(n)}$ together with all points in $\T$ that are of the form $r_j=\rho(s_i)$
for some $s_i$ in the partition $\pi^{(n)}$. Note that under this construction we have $r_0,r_1,...,r_m\in \T$. In fact, for any interval
$(s_{i-1}, s_i)$ in $(a)$, there is an identical interval $(r_{j-1}, r_j)$ in the partition $\sm^{(n)}$ corresponding to it;
for any interval $(s_{i-1}, s_i)$ in $(b)$, there are two intervals $(r_{j-2}, r_{j-1})$ and $(r_{j-1}, r_j)$ corresponding to it, so
that $r_{j-1}=\rho(s_i)$. And by \eqref{EstimateOfRemainderIntervals} we know that
$$\sum\li_{(s_{i-1}, s_i)\in (b), \, (r_{j-2}, r_{j-1}) \text{ is corresponding to it}}(r_{j-1}-r_{j-2})<\dfrac{n}{2^n} \ .$$
Note that the number $m$ depends on $n$ and the partition $\pi^{(n)}$: $m=m(n ,\pi^{(n)})$. In particular
$m\ra \infty$ as $n \ra \infty$. For
simplicity we will suppress this dependence later in our proof.

Let us recall the definition of deterministic and stochastic $\Dt$--integrals as defined in Section 2.
Let $\tilde{f}$ be the extension of $f$ that we have in \eqref{Extension}: for any
$t\in \T$, $$\tilde{f}(t,\om)=f(\sup[0,t]_\T, \om) \ .$$
Note that if $t\in \T$ is such that $\rho(t)=t$, then $\tilde{f}(t,\om)=f(t,\om)$,
otherwise if $t\in \T$ is such that $\rho(t)<t$, then $\tilde{f}(t,\om)=f(\rho(t), \om)$.
Thus we see that
$$\Prob\left(\lim\li_{n\ra \infty}\sum\li_{j=1}^{m} f(r_{j-1}, \om)(r_j-r_{j-1})=\play{\int_{t_1}^{t_2}\tilde{f}(s,\om)ds}\right)=1 \  ,$$
$$\Prob\left(\lim\li_{n\ra \infty}\sum\li_{j=1}^{m} f(r_{j-1}, \om)(W_{r_j}-W_{r_{j-1}})=\play{\int_{t_1}^{t_2}\tilde{f}(s,\om)dW_s}\right)=1 \  .$$

So it suffices to prove that
$$\Prob\left(\lim\li_{n\ra \infty}\left[\sum\li_{i=1}^n f(s_{i-1}, \om)(s_i-s_{i-1})
-\sum\li_{j=1}^m f(r_{j-1}, \om)(r_j-r_{j-1})\right]=0\right)=1$$
and
$$\Prob\left(\lim\li_{n\ra \infty}\left[
\sum\li_{i=1}^n f(s_{i-1}, \om)(W_{s_i}-W_{s_{i-1}})
-\sum\li_{j=1}^m f(r_{j-1}, \om)(W_{r_j}-W_{r_{j-1}})\right]=0
\right)=1 \  .$$

In fact, for any interval $(s_{i-1}, s_i)$ in class $(a)$, there exist an interval $(r_{j-1}, r_j)$ identical to the
interval $(s_{i-1}, s_i)$, so that
$$f(s_{i-1}, \om)(s_i-s_{i-1})-f(r_{j-1}, \om)(r_j-r_{j-1})=0$$
and
$$f(s_{i-1}, \om)(W_{s_i}-W_{s_{i-1}})-f(r_{j-1}, \om)(W_{r_j}-W_{r_{j-1}})=0 \ .$$
For any open interval $(s_{i-1}, s_i)$
in class $(b)$, there are two corresponding intervals $(r_{j-2}, r_{j-1})$ and $(r_{j-1}, r_j)$
 such that $r_{j-2}=s_{i-1}$, $r_{j-1}=\rho(s_i)$ and $r_j=s_i$. In this case

$$\begin{array}{ll}
&f(s_{i-1}, \om)(s_i-s_{i-1})-f(r_{j-1}, \om)(r_j-r_{j-1})-f(r_{j-2}, \om)(r_{j-1}-r_{j-2})
\\
=&f(s_{i-1}, \om)(s_i-s_{i-1})-f(\rho(s_i), \om)(s_i-\rho(s_i))-f(s_{i-1}, \om)(\rho(s_i)-s_{i-1})
\\
=&(f(s_{i-1}, \om)-f(\rho(s_i), \om))(s_i-\rho(s_i))
\end{array}$$
and
$$\begin{array}{ll}
&f(s_{i-1}, \om)(W_{s_i}-W_{s_{i-1}})-f(r_{j-1}, \om)(W_{r_j}-W_{r_{j-1}})-f(r_{j-2}, \om)(W_{r_{j-1}}-W_{r_{j-2}})
\\
=&f(s_{i-1}, \om)(W_{s_i}-W_{s_{i-1}})-f(\rho(s_i), \om)(W_{s_i}-W_{\rho(s_i)})-f(s_{i-1}, \om)(W_{\rho(s_i)}-W_{s_{i-1}})
\\
=&(f(s_{i-1}, \om)-f(\rho(s_i), \om))(W_{s_i}-W_{\rho(s_i)}) \ .
\end{array}$$

From the above calculations and the fact that we have
\eqref{EstimateOfRemainderIntervals} and $f$ is continuous on $[t_1, t_2]\cap \T$, together with the fact that
$s_{j-1}, \rho(s_j)\in \T$, $0\leq \rho(s_j)-s_{j-1}\leq \dfrac{1}{2^n}$, we see the claim follows.
$\square$

\

\textbf{Lemma 2.} (Convergence of quadratic variation of Brownian motion on time scale.)

\textit{Given a time scale $\T$ and $t_1, t_2\in \T$, $t_1<t_2$; a probability space
$(\Om, \cF, \Prob)$; a Brownian motion $\{W_t\}_{t\in \T}$ on the time scale $\T$.
Let any $L^2([t_1, t_2]_\T)$--progressively
 measurable random function $f(t, \om)$ on $\T$ be defined such that $\E f^2(t,\om)$ is uniformly
 bounded on $[t_1, t_2]$. Consider the families of partitions $\pi^{(n)}: t_1=s_0<s_1<...<s_n=t_2$,
$s_0, s_1, ...,s_n\in \T$, $\max\li_{i=1,2,...,n}(\rho(s_i)-s_{i-1})<\dfrac{1}{2^n}$. We classify all the intervals
$(s_{i-1}, s_i)$, $i=1,2,...,n$ into two classes $(a)$ and $(b)$ as before. Then we have }
\begin{equation}\label{ContinuousPartQuadraticVariation}
\Prob\left(\lim\li_{n\ra \infty}\left[\sum\li_{(a)} f(s_{i-1}, \om)(W_{s_{i}}-W_{s_{i-1}})^2
-\sum\li_{(a)}f(s_{i-1}, \om)(s_i-s_{i-1})\right]=0\right)=1 \ .
\end{equation}

\

\textbf{Proof.} We notice that for all intervals $(s_{i-1}, s_i)\in (a)$ we have $\rho(s_{i})=s_{i-1}$
and thus $s_i-s_{i-1}<\dfrac{1}{2^n}$. Let us denote that

$$V_n=\left[\sum\li_{(a)} f(s_{i-1}, \om)(W_{s_{i}}-W_{s_{i-1}})^2
-\sum\li_{(a)}f(s_{i-1}, \om)(s_i-s_{i-1})\right] \ .$$

Since $f(t,\om)$ is progressively measurable, we see that $f(s_{i-1}, \om)$ is independent of $W_{s_i}-W_{s_{i-1}}$. Thus

$$\begin{array}{ll}
\E V_n &= \play{\E \left[\sum\li_{(a)} f(s_{i-1}, \om)(W_{s_{i}}-W_{s_{i-1}})^2
-\sum\li_{(a)}f(s_{i-1}, \om)(s_i-s_{i-1})\right]}
\\
&\play{=\sum\li_{(a)} \E f(s_{i-1}, \om)(s_{i}-s_{i-1})-\sum\li_{(a)} \E f(s_{i-1}, \om)(s_i-s_{i-1})}
\\
&=0 \ .
\end{array}
$$

Furthermore

$$\begin{array}{ll}
\Var V_n &= \play{\E \left[\sum\li_{(a)} f(s_{i-1}, \om)(W_{s_{i}}-W_{s_{i-1}})^2
-\sum\li_{(a)}f(s_{i-1}, \om)(s_i-s_{i-1})\right]^2}
\\
&\play{=\E \sum\li_{(a)} f(s_{i-1}, \om)f(s_{j-1}, \om)[(W_{s_{i}}-W_{s_{i-1}})^2-(s_{i}-s_{i-1})]
\cdot[(W_{s_{j}}-W_{s_{j-1}})^2-(s_{j}-s_{j-1})]}
\\
&\play{=\sum\li_{(a)} \E f(s_{i-1}, \om)f(s_{j-1}, \om)[(W_{s_{i}}-W_{s_{i-1}})^2-(s_{i}-s_{i-1})]
\cdot[(W_{s_{j}}-W_{s_{j-1}})^2-(s_{j}-s_{j-1})]} \ .
\end{array}
$$

If $i<j$, then $f(s_{i-1}, \om)f(s_{j-1}, \om)[(W_{s_{i}}-W_{s_{i-1}})^2-(s_{i}-s_{i-1})]$ is independent of
$[(W_{s_{j}}-W_{s_{j-1}})^2-(s_{j}-s_{j-1})]$, so we have $\E f(s_{i-1}, \om)f(s_{j-1}, \om)[(W_{s_{i}}-W_{s_{i-1}})^2-(s_{i}-s_{i-1})]
\cdot[(W_{s_{j}}-W_{s_{j-1}})^2-(s_{j}-s_{j-1})]=0$. Similarly, for $i>j$ we also have
$\E f(s_{i-1}, \om)f(s_{j-1}, \om)[(W_{s_{i}}-W_{s_{i-1}})^2-(s_{i}-s_{i-1})]
\cdot[(W_{s_{j}}-W_{s_{j-1}})^2-(s_{j}-s_{j-1})]=0$. This implies that

$$\begin{array}{ll}
\Var V_n &=\play{\sum\li_{(a)} \E f^2(s_{i-1}, \om)[(W_{s_{i}}-W_{s_{i-1}})^2-(s_{i}-s_{i-1})]^2}
\\
&=\play{\sum\li_{(a)} \E f^2(s_{i-1}, \om) \E [(W_{s_{i}}-W_{s_{i-1}})^2-(s_{i}-s_{i-1})]^2}
\\
&=\play{\sum\li_{(a)} \E f^2(s_{i-1}, \om) \E [(W_{s_{i}}-W_{s_{i-1}})^4-2(W_{s_{i}}-W_{s_{i-1}})^2(s_{i}-s_{i-1})+(s_{i}-s_{i-1})^2]}
\\
&=\play{\sum\li_{(a)} \E f^2(s_{i-1}, \om) [3(s_{i}-s_{i-1})^2-2(s_{i}-s_{i-1})^2+(s_{i}-s_{i-1})^2]}
\\
&=\play{2\sum\li_{(a)} \E f^2(s_{i-1}, \om) (s_{i}-s_{i-1})^2}\leq \dfrac{1}{2^{n-1}} \left(\max\li_{s\in[t_1,t_2]}\E f^2(s,\om)\right) \sum\li_{(a)} (s_{i}-s_{i-1})\ra 0
\end{array}$$
as $n\ra \infty$. This together with the fact that $\E V_n=0$ for any $n$ implies the claim \eqref{ContinuousPartQuadraticVariation} of the Lemma. $\square$

\

The argument above leads us to an It\^{o}'s formula for $f(t, W_t)$. Making use of the same methods, one can derive a more general It\^{o}'s formula
for the solution $X_t$ to the $\Dt$--stochastic differential equation \eqref{SDeltaE}. We will not repeat the proof, but we will claim the
following Theorem.

\

\textbf{Theorem 3.}  (General It\^{o}'s formula) \textit{Let $X_t$ be the solution to the $\Dt$--stochastic
differential equation \eqref{SDeltaE}. Let any function $f: \T\times \R\ra \R$ be such that
$f^\Dt(t,x)$, $f^{\Dt^2}(t,x)$, $\dfrac{\pt f}{\pt x}(t,x)$, $\dfrac{\pt^2 f}{\pt x^2}(t,x)$,
$\dfrac{\pt f^\Dt}{\pt x}(t,x)$ and $\dfrac{\pt^2 f^\Dt}{\pt x^2}(t,x)$
are continuous on $\T\times \R$.
For any $t_1\leq t_2$, $t_1, t_2\in [0, \infty)_\T$ we have}

\begin{equation}\label{ItoFormulaGeneral}\begin{array}{ll}
&f(t_2,X_{t_2})-f(t_1, X_{t_1})
\\
=&\play{\int_{t_1}^{t_2}b(s, X_s)f^\Dt(s, X_s)\Dt s+
\int_{t_1}^{t_2}\sm(s, X_s)\dfrac{\pt f}{\pt x}(s,X_s)\Dt W_s
+\dfrac{1}{2}\int_{t_1}^{t_2}\sm^2(s, X_s)\dfrac{\pt^2 f}{\pt x^2}(s, X_s)\Dt s}
\\
&+\play{\sum\li_{I_k\in \cI, I_k\subset (t_1,t_2)} \left[f(s_{I_k}^+,W_{s_{I_k}^+})-f(s_{I_k}^+, W_{s_{I_k}^-})
\right.}
\\
&\play{\ \ \ \ \ \ \ \left.-\sm(s_{I_k}^-, W_{s_{I_k}^-})\dfrac{\pt f}{\pt x}(s_{I_k}^-,W_{s_{I_k}^-})(W_{s_{I_k}^+}-W_{s_{I_k}^-})
-\dfrac{1}{2}\sm^2(s_{I_k}^-, W_{s_{I_k}^-})\dfrac{\pt^2 f}{\pt x^2}(s_{I_k}^-,W_{s_{I_k}^-})(s_{I_k}^+-s_{I_k}^-)\right] \ .}
\end{array}\end{equation}

\

\section{The stochastic exponential on time scales.}

Our target in this section is to establish a closed--form formula for the \textit{stochastic exponential} in the case of general time scales
$\T$.

\

\textbf{Definition 3.} \textit{We say an adapted stochastic process $A(t)$ defined on the filtered probability space $(\Om, \cF_t, \Prob)$
is} stochastic regressive \textit{with respect to the Brownian motion $W_t$ on the time scale $\T$ if and only if for any right--scattered point
$t\in \T$ we have}
$$(1+A(t))(W_{\sm(t)}-W_t) \neq 0  \ , \ \text{ a.s.  \ for all } t\in \T  \ .$$
\textit{The set of stochastic regressive processes will be denoted by }$\cR_W$.

\

The following definition of a stochastic exponential was also introduced in \cite{[SD Exponential Bohner-Sanyal]}.

\

\textbf{Definition 4.} (Stochastic Exponential) \textit{Let $t_0\in \T$ and $A\in \cR_W$, then the unique solution
of the $\Dt$--stochastic differential equation}
\begin{equation}\label{SDeltaEStochasticExponential}
\Dt X_t=A(t)X_t\Dt W_t \ , \ X(t_0)=1 \ , \ t\in \T
\end{equation}
\textit{is called the} stochastic exponential \textit{and is denoted by}
$$X_\bullet=\cE_A(\bullet, t_0) \  .$$

\

We note that $\cE_A(t, t_0)$ as a solution to the equation
\eqref{SDeltaEStochasticExponential} can be written into an integral equation

\begin{equation}\label{StochasticExponentialIntagralEquation}
\cE_A(t, t_0)=1+\int_{t_0}^t A(s)\cE_A(s,t_0)\Dt W_s \ , \ \text{ for all } t\in \T \  .
\end{equation}
We will be making use of the set--up we have in Section 4 about It\^{o}'s formula.
Let $t_0<t$ and $t_0, t\in \T$.
Let the sets $C$ and $\cI$ be defined as in Section 4 corresponding to the interval $[t_1, t_2]=[t_0, t]$.
Let $I_k\in \cI$ and $I_k=(s_{I_k}^-, s_{I_k}^+)$. We note that $s_{I_k}^-=\rho(s_{I_k}^+)$, $s_{I_k}^+=\sm(s_{I_k}^-)$. Let
\begin{equation}\label{D}
D(t, t_0)=\sum\li_{I_k\in \cI, I_k\subset (t_0,t)}A(s_{I_k}^-)(W_{s_{I_k}^+}-W_{s_{I_k}^-})
-\dfrac{1}{2}\sum\li_{I_k\in \cI, I_k\subset (t_0,t)}A^2(s_{I_k}^-)(s_{I_k}^+-s_{I_k}^-) \ .
\end{equation}

We define
\begin{equation}\label{U}
U(t, t_0)=\prod\limits_{I_k\in \cI, I_k\subset (t_0,t)}\left[1+A(s_{I_k}^-)(W_{s_{I_k}^+}-W_{s_{I_k}^-})\right] \ ,
\end{equation}

\begin{equation}\label{V}
V(t, t_0)=\exp\left(\play{\int_{t_0}^t A(s)\Dt W_s-\dfrac{1}{2}\int_{t_0}^t A^2(s)\Dt s
- D(t, t_0)}\right) \ .
\end{equation}

\

\textbf{Theorem 4.} (Stochastic Exponential on time scales)
\textit{The stochastic exponential has the closed--form expression}
$$\cE_A(t, t_0)=U(t, t_0)V(t, t_0) \ .$$

\

\textbf{Proof.} Consider the process
$$Y_t=\int_{t_0}^t A(s)\Dt W_s-\dfrac{1}{2}\int_{t_0}^t A^2(s)\Dt s
- D(t, t_0) \ .$$

Let us introduce another function $\al(t)$ such that
$$\al(t)=\left\{
\begin{array}{ll}
0 & \ , \ \text{when } t=s_{I_k}^- \text{ or } t=s_{I_k}^+ \ ,
\\
A(t) & \ , \ \text{otherwise} \ .
\end{array}\right.$$

We see now that the process $Y_t$ is a solution to the $\Dt$--stochastic differential equation
$$\Dt Y_t=\al(t)\Dt W_s-\dfrac{1}{2}\al^2(t)\Dt s \ , \ Y_{t_0}=0 \ .$$

Notice that $Y_{s_{I_k}^-}=Y_{s_{I_k}^+}$ for any $I_k=(s_{I_k}^-, s_{I_k}^+)\in \cI$. Taking this into account, as well
as the fact that $\al(s)=0$ whenever $t=s_{I_k}^- \text{ or } t=s_{I_k}^+$,
we can apply the general It\^{o}'s formula \eqref{ItoFormulaGeneral} to the function $V(t, t_0)=\exp(Y_t)$
and we will get

$$\begin{array}{ll}
&\exp(Y_t)-1
\\
=&\play{\int_{t_0}^t \al(s)\exp(Y_s)\Dt W_s-\dfrac{1}{2}\int_{t_0}^t \al^2(s)\exp(Y_s)\Dt s+\dfrac{1}{2}\int_{t_0}^t\al^2(s)\exp(Y_s)\Dt s}
\\
=&\play{\int_{t_0}^t \al(s)\exp(Y_s)\Dt W_s} \ .
\end{array}$$

Thus
$$V(t, t_0)=1+\int_{t_0}^t \al(s)V(s, t_0)\Dt W_s \ ,$$
or in other words
$$\Dt V(t, t_0)=\al(t)V(t, t_0)\Dt W_t \ .$$

Let us now consider the function $\cE_A(t, t_0)=U(t, t_0)V(t, t_0)$. We claim that

\begin{equation}\label{Claim}
U(t, t_0)V(t, t_0)-1=\int_{t_0}^t A(s)U(s,t_0)V(s, t_0)\Dt W_s \ .
\end{equation}

Notice that
$$U(s_{I_k}^+, t_0)=[1+A(s_{I_k}^-)(W_{s_{I_k}^+}-W_{s_{I_k}^-})]U(s_{I_k}^-, t_0)=
U(s_{I_k}^-, t_0)+A(s_{I_k}^-)U(s_{I_k}^-, t_0)(W_{s_{I_k}^+}-W_{s_{I_k}^-}) \ ,$$
i.e.,
$$U(s_{I_k}^+, t_0)-U(s_{I_k}^-, t_0)=
A(s_{I_k}^-)U(s_{I_k}^-, t_0)(W_{s_{I_k}^+}-W_{s_{I_k}^-}) \ .$$
Using this fact,
the above claimed identity \eqref{Claim} can be written as

\begin{equation}\label{ProductDifferenceToProve}
\begin{array}{ll}
&U(t, t_0)V(t, t_0)-1
\\
=&\play{\int_{t_0}^t \al(s)U(s,t_0)V(s,t_0)\Dt W_s
+\sum\li_{I_k\in \cI, I_k\subset (t_0, t)}A(s_{I_k}^-)U(s_{I_k}^-,t_0)V(s_{I_k}^-,t_0)(W_{s_{I_k}^+}-W_{s_{I_k}^-})}
\\
=&\play{\int_{t_0}^t \al(s)U(s,t_0)V(s,t_0)\Dt W_s
+\sum\li_{I_k\in \cI, I_k\subset (t_0, t)}V(s_{I_k}^-,t_0)(U(s_{I_k}^+, t_0)-U(s_{I_k}^-, t_0))} \ .
\end{array}
\end{equation}

In fact, with respect to the partition $\pi^{(n)}: t_0=s_0<s_1<...<s_{n-1}<s_n=t$ that we have been using, we have

$$\begin{array}{ll}
&U(t, t_0)V(t, t_0)-1
\\
=&\play{\sum\li_{i=1}^n [U(s_i, t_0)V(s_i, t_0)-U(s_{i-1}, t_0)V(s_{i-1}, t_0)]}
\\
=&\play{\sum\li_{i=1}^n \left[(U(s_i, t_0)-U(s_{i-1}, t_0))(V(s_i, t_0)-V(s_{i-1}, t_0))
+U(s_{i-1}, t_0)(V(s_i, t_0)-V(s_{i-1}, t_0))\right.}
\\
& \ \ \ \ \ \ \ \ \ \ \ \ \ \ \ \play{\left.+V(s_{i-1}, t_0)(U(s_i, t_0)-U(s_{i-1}, t_0))\right]}
\\
=&\play{\sum\li_{i=1}^n (U(s_i, t_0)-U(s_{i-1}, t_0))(V(s_i, t_0)-V(s_{i-1}, t_0))
+\sum\li_{i=1}^n U(s_{i-1}, t_0)(V(s_i, t_0)-V(s_{i-1}, t_0))}
\\
& \ \ \ \ \ \ \ \play{+\sum\li_{i=1}^n V(s_{i-1}, t_0)(U(s_i, t_0)-U(s_{i-1}, t_0))}
\\
=&(I)+(II)+(III) \ .
\end{array}$$

Here

$$\begin{array}{llll}
(I)&=& \play{\sum\li_{i=1}^n (U(s_i, t_0)-U(s_{i-1}, t_0))(V(s_i, t_0)-V(s_{i-1}, t_0))} & \ ,
\\
(II)&=& \play{\sum\li_{i=1}^n U(s_{i-1}, t_0)(V(s_i, t_0)-V(s_{i-1}, t_0))} & \ ,
\\
(III)&=& \play{\sum\li_{i=1}^n V(s_{i-1}, t_0)(U(s_i, t_0)-U(s_{i-1}, t_0))} & \ .
\end{array}$$

We can apply the previous arguments and classify the intervals $(s_{i-1}, s_i)$ into classes $(a)$
and $(b)$. Notice that on each interval $(s_{I_k}^-, s_{I_k}^+)$, the function $V(t, t_0)$ remains constant and the function
$U(t, t_0)$ has a jump, and on each interval $(s_{i-1}, s_i)$ in class $(a)$ the function $U(t, t_0)$ is a constant. This observation
and similar arguments (which we leave to the reader) as in the previous section will enable us to prove that with probability $1$, as $n \ra \infty$,
we will have

$$\begin{array}{llll}
(I)  &\ra &0 & \ , \\
(II) &\ra &\play{\int_{t_0}^t \al(s)U(s,t_0)V(s,t_0)\Dt W_s} & \ , \\
(III)&\ra &\play{\sum\li_{I_k\in \cI, I_k\subset (t_0, t)}V(s_{I_k}^-,t_0)(U(s_{I_k}^+, t_0)-U(s_{I_k}^-, t_0))} & \ .
\end{array}$$

 So we proved
\eqref{ProductDifferenceToProve} and thus \eqref{Claim}. $\square$

\section{Change of measure and Girsanov's theorem on time scales.}

We demonstrate in this section
a change of measure formula (Girsanov's formula) for Brownian motion on time scales. Our analysis is based
on the method of extension that was introduced in Section 3 (originally from \cite{[SDtE Bohner]}).

Let us consider two processes, the standard Brownian motion $\{W_t\}_{t\in \T}$ on $(\Om, \cF_t, \Prob)$
on the time scale $\T$, and the process

$$B_t=W_t-\int_0^t A(s)\Dt s \ ,$$
on the time scale $t\in \T$.

Let us consider an extension of the (probably random) function $A(s)$ as in \eqref{Extension}. Let us define the
so obtained extension function to be $\widetilde{A}(s)$. Recall that \eqref{Extension} implies that
$$\widetilde{A}(s,\om)=A(\sup[0,s]_\T,\om) \ .$$
Let $\widetilde{W}_t$ be a standard Brownian motion on $[0,\infty)$. If we define
$$\widetilde{B}_t=\widetilde{W}_t-\int_0^t \widetilde{A}(s)ds \ ,$$
then the process $\widetilde{B}_t$ agrees with $B_t$ for any time point $t\in \T$.

For any $t, t_0\in \T$, $t>t_0$, let

\begin{equation}\label{GirsanovDensity}
\begin{array}{ll}
\cG_A(t, t_0)&=\play{\exp\left(\int_{t_0}^t \widetilde{A}(s)d W_s-\dfrac{1}{2}\int_{t_0}^t \widetilde{A}^2(s)d s\right)}
\\
&=\play{\exp\left(\int_{t_0}^t \widetilde{A}(s)d W_s-\dfrac{1}{2}\int_{t_0}^t \widetilde{A^2}(s)d s\right)}
\\
&=\play{\exp\left(\int_{t_0}^t A(s)\Dt W_s-\dfrac{1}{2}\int_{t_0}^t A^2(s)\Dt s\right) \ .}
\end{array}
\end{equation}

It is easy to see that the function $\cG_A(t, t_0)$ is the standard Girsanov's density function for the process $\widetilde{B}_t$ with respect to
the standard Brownian motion $\widetilde{W}_t$. Since $\widetilde{B}_t$ and $\widetilde{W}_t$ have the same distributions
as $B_t$ and $W_t$ on the time scale $\T$, we conclude with the following two Theorems.

\

\textbf{Theorem 5.} (Novikov's condition on time scales)
\textit{If for every $t\geq 0$ we have }
\begin{equation}\label{NovikovCondition}
\E \exp\left(\int_0^t A^2(s)\Dt s\right)<\infty \ ,
\end{equation}
\textit{then for every $t\geq 0$ we have }
$$\E \cG_A(t, t_0)=1 \ .$$

\

Let \eqref{NovikovCondition} be satisfied. Let $T>0$ and pick $T>t_0$, $t_0, T\in \T$. Consider a new measure $\Prob^B$ on $(\Om, \cF_t)$,
defined by it Radon--Nikodym derivative with respect to $\Prob^W$, as
$$\dfrac{d\Prob^B}{d\Prob^W}=\cG_A(T, t_0) \ .$$

\

\textbf{Theorem 6.} (Girsanov's change of measure on time scales)
\textit{Under the measure $\Prob^B$ the process $B_t$, $t\in [0,T]_\T$
is a standard Brownian motion on $\T$.}

\

\section{Application to Brownian motion on a quantum time scale.}

In this section we are going to apply our result to a quantum time scale ($q$--time scale, see
\cite[Example 1.41]{[Bohner Peterson book]}). Let $q>1$ and
$$q^\Z:=\{q^k: k\in \Z\} \text{ and } \overline{q^\Z}:=q^\Z\cup \{0\} \ .$$

The quantum time scale ($q$--time scale) is defined by $\T=\overline{q^\Z}$. Given the quantum time scale $\T$,
one can then construct a Brownian motion $W_t$ on $\T$ according to Definition 2.

We have
$$\sm(t)=\inf\{q^n: n\in [m+1,\infty)\}=q^{m+1}=qq^m=qt$$
if $t=q^m\in \T$ and obviously $\sm(0)=0$. So we obtain
$$\sm(t)=qt \text{ and } \rho(t)=\dfrac{t}{q} \text{ for all } t\in \T$$
and consequently $$\mu(t)=\sm(t)-t=(q-1)t \ \text{ for all } t\in \T\ .$$

Here $0$ is a right--dense minimum and every other point in $\T$ is isolated. For a function $f: \T\ra \R$
we have

$$f^\Dt(t)=\dfrac{f(\sm(t))-f(t)}{\mu(t)}=\dfrac{f(qt)-f(t)}{(q-1)t} \text{ for all } t\in \T\backslash \{0\}$$
and
$$f^\Dt(0)=\lim\li_{s\ra 0}\dfrac{f(0)-f(s)}{0-s}=\lim\li_{s\ra 0}\dfrac{f(s)-f(0)}{s}$$
provided the limit exists.

The open intervals $I_k$ that we have constructed in Section 4 have the form $I_k=(q^k, q^{k+1})$
where $k\in \Z$. For any two points $t_1<t_2$, $t_1,t_2 \in \T$, if $t_1, t_2\neq 0$, then $t_1=q^{k_1}$
and $t_2=q^{k_2}$ for two integers $k_1<k_2$. In this case we can apply \eqref{ItoFormula} and we get

\begin{equation}\label{ItoFormula1BMqScale}
\begin{array}{ll}
&f(q^{k_2}, W_{q^{k_2}})-f(q^{k_1}, W_{q^{k_1}})
\\
=&\play{\int_{q^{k_1}}^{q^{k_2}}f^\Dt(s, W_s)\Dt s+
\int_{q^{k_1}}^{q^{k_2}}\dfrac{\pt f}{\pt x}(s,W_s)\Dt W_s
+\dfrac{1}{2}\int_{q^{k_1}}^{q^{k_2}}\dfrac{\pt^2 f}{\pt x^2}(s, W_s)\Dt s}
\\
&\qquad +\play{\sum\li_{k=k_1}^{k_2-1} \left[f(q^{k+1},W_{q^{k+1}})-f(q^{k+1}, W_{q^k})
\right.}
\\
&\play{\qquad \qquad \qquad \qquad \qquad \left.-\dfrac{\pt f}{\pt x}(q^k,W_{q^k})(W_{q^{k+1}}-W_{q^k})
-\dfrac{1}{2}\dfrac{\pt^2 f}{\pt x^2}(q^k,W_{q^{k}})(q^{k+1}-q^k)\right] \ .}
\end{array}
\end{equation}

Since $\T\backslash\{0\}$ is a discrete time scale, we have

$$\int_{q^{k_1}}^{q^{k_2}}\dfrac{\pt f}{\pt x}(s,W_s)\Dt W_s=
\sum\li_{k=k_1}^{k_2-1} \dfrac{\pt f}{\pt x}(q^k,W_{q^k})(W_{q^{k+1}}-W_{q^k}) \ ,$$
and
$$\dfrac{1}{2}\int_{q^{k_1}}^{q^{k_2}}\dfrac{\pt^2 f}{\pt x^2}(s, W_s)\Dt s
=\sum\li_{k=k_1}^{k_2-1}\dfrac{1}{2}\dfrac{\pt^2 f}{\pt x^2}(q^k,W_{q^{k}})(q^{k+1}-q^k) \ .$$

Moreover, we have

$$\begin{array}{ll}
&\play{\int_{q^{k_1}}^{q^{k_2}}f^\Dt(s, W_s)\Dt s}
\\
=&\play{\sum\li_{k=k_1}^{k_2-1}\dfrac{f(q^{k+1},W_{q^k})-f(q^k, W_{q^k})}{q^{k+1}-q^k}(q^{k+1}-q^{k})}
\\
=&\play{\sum\li_{k=k_1}^{k_2-1}[f(q^{k+1},W_{q^k})-f(q^k, W_{q^k})] \ .}
\end{array}$$

Therefore \eqref{ItoFormula1BMqScale} becomes
$$
\begin{array}{ll}
&f(q^{k_2}, W_{q^{k_2}})-f(q^{k_1}, W_{q^{k_1}})
\\
=&\play{\sum\li_{k=k_1}^{k_2-1} \left\{[f(q^{k+1},W_{q^k})-f(q^k, W_{q^k})]+[f(q^{k+1},W_{q^{k+1}})-f(q^{k+1}, W_{q^k})]\right\}}
\\
=&\play{\sum\li_{k=k_1}^{k_2-1} [f(q^{k+1},W_{q^{k+1}})-f(q^k, W_{q^k})] \ ,}
\end{array}
$$
which is a trivial telescoping identity. This justifies \eqref{ItoFormula} in the case away from $0$.

Let us consider now the case when $t_1=0$ and $t_2=q^k>0$ for some $k\in \Z$. In this case we have, according to \eqref{ItoFormula},
that

\begin{equation}\label{ItoFormula2BMqScale}
\begin{array}{ll}
&f(q^{k}, W_{q^{k}})-f(0,0)
\\
=&\play{\int_{0}^{q^{k}}f^\Dt(s, W_s)\Dt s+
\int_{0}^{q^{k}}\dfrac{\pt f}{\pt x}(s,W_s)\Dt W_s
+\dfrac{1}{2}\int_{0}^{q^{k}}\dfrac{\pt^2 f}{\pt x^2}(s, W_s)\Dt s}
\\
&\qquad +\play{\sum\li_{j=-\infty}^{k-1} \left[f(q^{j+1},W_{q^{j+1}})-f(q^{j+1}, W_{q^j})
\right.}
\\
&\play{\qquad \qquad \qquad \qquad \qquad \left.-\dfrac{\pt f}{\pt x}(q^j,W_{q^j})(W_{q^{j+1}}-W_{q^j})
-\dfrac{1}{2}\dfrac{\pt^2 f}{\pt x^2}(q^j,W_{q^{j}})(q^{j+1}-q^j)\right] \ .}
\end{array}
\end{equation}

One can justify that in this case we have

$$\int_{0}^{q^{k}}\dfrac{\pt f}{\pt x}(s,W_s)\Dt W_s=
\sum\li_{j=-\infty}^{k-1} \dfrac{\pt f}{\pt x}(q^j,W_{q^j})(W_{q^{j+1}}-W_{q^j}) \ ,$$
and
$$\dfrac{1}{2}\int_{0}^{q^{k}}\dfrac{\pt^2 f}{\pt x^2}(s, W_s)\Dt s
=\sum\li_{j=-\infty}^{k-1}\dfrac{1}{2}\dfrac{\pt^2 f}{\pt x^2}(q^j,W_{q^{j}})(q^{j+1}-q^j) \ .$$

Moreover, we have

$$\begin{array}{ll}
&\play{\int_{0}^{q^{k}}f^\Dt(s, W_s)\Dt s}
\\
=&\play{\sum\li_{j=-\infty}^{k-1}\dfrac{f(q^{j+1},W_{q^j})-f(q^j, W_{q^j})}{q^{j+1}-q^j}(q^{j+1}-q^{j})}
\\
=&\play{\sum\li_{j=-\infty}^{k-1}[f(q^{j+1},W_{q^j})-f(q^j, W_{q^j})] \ .}
\end{array}$$

Therefore \eqref{ItoFormula1BMqScale} becomes
$$
\begin{array}{ll}
&f(q^{k}, W_{q^{k}})-f(0, 0)
\\
=&\play{\sum\li_{j=-\infty}^{k-1} \left\{[f(q^{j+1},W_{q^j})-f(q^j, W_{q^j})]+[f(q^{j+1},W_{q^{j+1}})-f(q^{j+1}, W_{q^j})]\right\}}
\\
=&\play{\sum\li_{j=-\infty}^{k-1} [f(q^{j+1},W_{q^{j+1}})-f(q^j, W_{q^j})] \ .}
\end{array}
$$
which is also a telescoping identity. This justifies \eqref{ItoFormula} in the case including $0$.

Making use of Theorem 4, it is easy to write down the stochastic exponential for the quantum time scale:

$$\cE_A(0,q^k)=\prod\limits_{j=-\infty}^{k-1}\left[1+A(q^j)(W_{q^{j+1}}-W_{q^j})\right] \ .$$

\

\textbf{Acknowledgements}:
The author would like to express his sincere gratitude to Professor Martin Bohner
for inviting him to present this work at the Time Scales Seminar, 
Missouri S\&T on October 5, 2016, and also for 
pointing out the use of delta (Hilger) derivatives. 
He would also like to thank Professor David Grow
for an inspiring discussion on the topic. Many thanks go to the anonymous referee for invaluable comments
that lead to essential improvements of the paper.
Special thanks are dedicated to Missouri University of Science and Technology (formerly University of Missouri, Rolla)
for a startup fund that supports this work.

\end{document}